\documentclass[12pt, final]{amsart}
\usepackage[english]{babel}
\usepackage{latexsym}
\usepackage{amssymb}
\usepackage{amscd}
\usepackage[cp850]{inputenc}
\usepackage[mathscr]{eucal}
\title{On Lindenstrauss-Pe\l czy\'{n}ski spaces}

\author{Jes\'us M. F. Castillo, Yolanda Moreno and Jes\'us Su\'arez}

\address{Departamento de Matemáticas, Universidad de
Extremadura, Avenida de Elvas, 06071 Badajoz, Espana. e-mail
address: castillo@unex.es; ymoreno@unex.es,
jesussf@telefonica.net}

\thanks{The research of the first two authors has been supported in part by DGICYT project MTM2004-02635. The work of the
third author was supported in part by a Marie Curie grant
HPMT-GH-01-00286-04 at Karlsruhe University under the direction of
Prof. L. Weis.}


\date{}

\begin{document}
\maketitle \markboth{}{} \theoremstyle{plain}
\newtheorem{theorem}{Theorem}
\newtheorem{defin}{Definition}
\newtheorem{proposition}{Proposition}[section]
\newtheorem{prop}{Proposition}[section]
\newtheorem{corollary}{Corollary}
\newtheorem{cor}{Corollary}
\newtheorem{lemma}{Lemma}
\newtheorem{lema}{Lemma}
\theoremstyle{remark}
\newtheorem{problem}{Problem}
\newcommand{\R}{\mathbb{R}}
\newcommand{\N}{\mathbb{N}}
\newcommand{\seq}{\ensuremath{0\to Y\to X\to Z\to0\:}}
\newcommand{\seqq}{\ensuremath{0\to Y_1\to X_1\to Z_1\to0\:}}
\newcommand{\seqF}{\ensuremath{0\to Y\to Y\oplus_FZ\to Z\to0\:}}
\newcommand{\seqG}{\ensuremath{0\to Y\to Y\oplus_GZ\to Z\to0\:}}
\newcommand{\F}{\ensuremath{F:Z\to Y\:}}
\newcommand{\Lin}{\ensuremath{L:Z\to Y\:}}
\newcommand{\B}{\ensuremath{B:Z\to Y\:}}
\newcommand{\YFZ}{\ensuremath{Y\oplus_FZ\:}}

\def\Ext{\operatorname{Ext}}
\bibliographystyle{plain}
\begin{abstract} In this work we shall be concerned with some stability aspects of
the classical problem of extension of $C(K)$-valued operators. We
introduce the class $\mathscr{LP}$ of Banach spaces of
Lindenstrauss-Pe\l czy\'{n}ski type as those such that every
operator from a subspace of $c_0$ into them can be extended to
$c_0$. We show that all $\mathscr{LP}$-spaces are of type
$\mathcal L_\infty$ but not the converse. Moreover, $\mathcal
L_\infty$-spaces will be characterized as those spaces $E$ such
that $E$-valued operators from $w^*(l_1,c_0)$-closed subspaces of
$l_1$ extend to $l_1$. Complemented subspaces of $C(K)$ and
separably injective spaces are subclasses of $\mathscr{LP}$-spaces
and we show that the former does not contain the latter. It is
established that $\mathcal L_\infty$-spaces not containing $l_1$
are quotients of $\mathscr{LP}$-spaces, while $\mathcal
L_\infty$-spaces not containing $c_0$, quotients of an
$\mathscr{LP}$-space by a separably injective space and twisted
sums of $\mathscr{LP}$-spaces are $\mathscr{LP}$-spaces.
\end{abstract}
\section{Introduction and Preliminaries}
In this work we shall be concerned with some stability aspects of
the classical problem of extension of $C(K)$-valued and $\mathcal
L_\infty$-valued operators. Let us describe and motivate them. In
a 1971 paper \cite{lindpelc} Lindenstrauss and Pe\l czy\'{n}ski
proved:
\begin{theorem} Let $K$ be a compact Hausdorff space.
Every $C(K)$-valued operator defined on a subspace of $c_0$ admits
an extension to the whole space.
\end{theorem}
The result remained isolated for quite a long time until 1989 when
Johnson and Zippin obtained in \cite{johnzipproc} an extension to
subspaces of $c_0(\Gamma)$, and later in 1995, in \cite{jz*}, the
analogous result for $w(l_1,c_0)$-closed subspaces of $l_1$.
Further proofs of the Lindenstrauss-Pe\l czy\'{n}ski theorem have
been provided by Zippin \cite{zippln,zipproc}. The paper
\cite{ccmy} contains an homological approach to both results
showing that they are in a sense dual one of the other.\smallskip

The general problem of extension of operators admits a natural
formulation in homological terms. We shall assume from the reader
some familiarity with the basic notions and constructions of the
theory of exact sequences of Banach spaces; the necessary
background can be seen in \cite{castgonz} and, operatively
defined, below. We shall write $0\to Y\stackrel{j}\to X
\stackrel{q}\to Z\to 0\equiv F$ to represent an exact sequence of
Banach spaces and operators, which is a diagram where the kernel
of each operator coincides with the image of the preceding. The
open mapping theorem makes $Y$ a subspace of $X$ through the
embedding $j$ and $Z$ the corresponding quotient space through
$q$. The reader can understand that $F$ is just the name of the
sequence, which is everything one needs to follow this paper;
however, those familiar with the theory of quasi-linear maps
created in \cite{kalt,kaltpeck} can in fact consider $F$ as a
quasi-linear map associated to the exact sequence.\smallskip

We shall consider exact sequences of Banach spaces modulus the
natural equivalence relationship: two sequences $F$ and $G$ are
said to be equivalent if there is a commutative diagram$$
\begin{CD}0@>>> Y@>j>> X@>>> Z@>>>0\equiv F\\
&& \Vert&& @VTVV \Vert\\
0@>>> Y@>>> X'@>>> Z@>>>0\equiv G.
\end{CD}
$$ In this case we write $F\equiv G$. The space of equivalence
classes of exact sequences with Y as subspace and $Z$ as quotient
will be denoted $\Ext(Z,Y)$. It is a vector space under some
natural operations (see \cite{mac}), and the $0$ element is the
sequence $0\to Y\to Y\oplus Z\to Z\to 0$ with inclusion $y\to
(y,0)$ and quotient map $(y,z)\to z$. We shall say that $F$ is
trivial or splits when $F\equiv 0$. This means, in classical
terms, that $j(Y)$ is complemented in $X$. Recall that a property
$P$ is said to be a 3-space property if whenever one has an exact
sequence  $0\to Y \to X \to Z \to 0$ in which both $Y$ and $Z$
have $P$ then also $X$ has $P$; see \cite{castgonz} for general
information about 3-space problems.

The lower sequence in a diagram
$$\begin{CD}0@>>> Y@>j>> X@>q>> Z@>>>0\equiv F\\ &&\|&& @AAA @AATA
\\ 0@>>> Y@>>> X'@>>> E@>>>0,\qquad
\end{CD}
$$is called {\it pull-back sequence} and naturally denoted $FT$.
The sequence $FT$ splits if and only if $T$ can be lifted to $X$
through $q$. The lower sequence in a diagram
$$
\begin{CD}0@>>> Y@>j>> X@>>> Z@>>>0\equiv F\\
&& @VTVV@VVV \Vert\\
0@>>> E@>>> X'@>>> Z@>>>0
\end{CD}
$$is called the {\it push-out sequence} and naturally denoted
$TF$. Extending an operator $T:Y\to E$ through $j$ is the same as
saying that $T F$ is trivial. The lifting property of
$l_1(\Gamma)$ and the fact that every Banach space $Z$ admits an
exact sequence $0 \to K(Z) \to l_1(\Gamma) \to Z \to 0$, called
{\it projective presentation} of $Z$, yield that every exact
sequence $0\to Y \to X \to Z\to 0$ is a push-out of a projective
presentation of $Z$. Hence $\Ext(Z,Y)=0$ is equivalent to the
statement ``every operator $K(Z) \to Y$ can be extended to
$l_1(\Gamma)$''.\smallskip

That all operators $Y\to E$ can be extended to $X$ through $j$
admits an even simpler formulation: the restriction operator
$j^*:\mathcal L(X,E)\to \mathcal L(Y,E)$ is surjective. The
following notation is quite natural and will prove to be very
useful:
\begin{defin}[$\mathcal A$-trivial exact sequences] Let $\mathcal A$ be a class of Banach spaces.
We say that an exact sequence $0\to Y\stackrel{j}\to X\to Z\to
0\equiv F$ is $\mathcal A$-trivial (or that $\mathcal A$-splits)
if for every $A \in \mathcal A$ the restriction operator
$j^*:\mathcal L(X, A)\to \mathcal L(Y, A)$ is
surjective.\end{defin} Sometimes the quantitative version of the
previous notion shall be required: given $\lambda\geq 1$, the
exact sequence $F$ will be said to be $(\lambda, \mathcal
A)$-trivial if for every $A \in \mathcal A$ every operator $T:Y\to
A$ admits an extension $\widehat{T}:X\to A$ such that
$\|\widehat{T}\|\leq \lambda \|T\|$.\smallskip

This notion of $\mathcal A$-triviality unifies different notions
appearing in the literature: i) trivial sequences, which
correspond with $\mathcal A$ = all Banach spaces; ii) Kalton's
{\it locally trivial}, or locally split, sequences (see
\cite{kaltloc}), corresponding with $\mathcal A = l_\infty(G_n)$,
where $G_n$ is a dense (in the Banach-Mazur distance) sequence of
finite dimensional Banach spaces (see also \cite{john}); iii)
Zippin's {\it almost trivial} sequences (see
\cite{zippill,zippln,zipproc,zipphandbook}), which correspond with
the choice $\mathcal A$= $C(K)$-spaces.\smallskip

In this work we are concerned with $C(K)$-trivial and $\mathcal
L_\infty$-trivial sequences. In Section 2 we study the stability
of $C(K)$-trivial sequences by amalgams and duality. We first show
that $l_p$ and $c_0$-amalgams of $C(K)$-trivial sequences are
$C(K)$-trivial. Regarding stability of $C(K)$-trivial sequences by
duality, the Lindenstrauss-Pe\l czy\'{n}ski and Johnson-Zippin
theorems suggest that it could be that the dual of a
$C(K)$-trivial sequence is $C(K)$-trivial since the former implies
that every exact sequence $0\to H \to X \to S \to 0 \equiv F$ in
which $H$ is a subspace of $c_0$ and $S$ separable is
$C(K)$-trivial; and the latter yields (see \cite{castmoreisr})
that its dual sequence $F^*$ is $C(K)$-trivial. However, the
situation outlined by those two theorems proves to be quite
peculiar; we give examples at the end of Section 2 which show that
the dual and bidual sequences of a $C(K)$-trivial sequence need
not to be $C(K)$-trivial.

In Section 3 our attention turns to those Banach spaces which can
play the role of $C(K)$-spaces in the Lindenstrauss-Pe\l
czy\'{n}ski theorem. We give to such spaces the name of
Lindenstrauss-Pe\l czy\'{n}ski ($\mathscr{LP}$, in short) spaces.
Our motivation to introduce those spaces comes from \cite[remark
2, p.234]{lindpelc} in which Lindenstrauss and Pe\l czy\'{n}ski
assert that isometric $L_1$-preduals can play the role of
$C(K)$-spaces regarding extension of operators from subspaces of
$c_0$. After showing that every $\mathscr{LP}$-space is an
$\mathcal L_\infty$ space we then face the unavoidable question:
Must every $\mathcal L_\infty$-space be an $\mathscr{LP}$-space?
The answer is no, which solves Problem 6.15 of Zippin in
\cite{zipphandbook}. We shall show that the same approach with
respect to the Johnson-Zippin theorem just provides (see Prop.
\ref{jz}) a new characterization of $\mathcal L_\infty$-spaces.

Returning to the problem of identification of $\mathscr{LP}$
spaces, it is clear that complemented subspaces of $C(K)$-spaces
and separably injective spaces are $\mathscr{LP}$-spaces. By the
way, we will show which is perhaps the first example of a
separably injective space not complemented in any $C(K)$-space.
The previous examples do not exhaust the class of
$\mathscr{LP}$-spaces: we shall show that $\mathcal
L_\infty$-spaces not containing $c_0$, the new exotic $\mathcal
L_\infty$-spaces constructed in \cite{ccky}, the quotients of
$\mathscr{LP}$-spaces by separably injective subspaces and the
$c_0$-vector sums of {\it uniformly} $\mathscr{LP}$-spaces are
$\mathscr{LP}$-spaces.

In Section 4 we tackle the 3-space problem for the class of
$\mathscr{LP}$-spaces, which needs the development of a new method
of proof and new characterizations of $\mathscr{LP}$-spaces.
Section 5 contains frther remarks, examples and open problems.

\section{On the stability of $C(K)$-trivial sequences by amalgams and duality}
To carry on our study on the stability of $C(K)$-trivial sequences
we need to know their behavior with respect to the basic
homological pull-back and push-out constructions.
\begin{proposition}\label{pbct} Let $\mathcal A$ be a class of
Banach spaces.
\begin{enumerate}
\item A pull-back sequence of an $\mathcal A$-trivial sequence is
$\mathcal A$-trivial. \item If $F$ is an $\mathcal A$-trivial
sequence and $\phi$ is a surjective operator then the push-out
sequence $\phi F$ is $\mathcal A$-trivial.
\end{enumerate}\end{proposition}
\begin{proof}The first assertion is obvious. The second is consequence of the next lemma.
\end{proof}\enlargethispage*{50pt}
\begin{lemma}\label{cruz}Let $\mathcal A$ be a class of Banach spaces.
Consider the completed push-out diagram of Banach spaces
$$
\begin{CD}
 & &0 &&0\\
 & &@VVV @VVV\\
 & &B&=& B\\
 & &@VV{a}V  @VV{b}V\\
 0@>>> Y @>j>> X @>>> Z @>>> 0\equiv F\\
 & &@VV{c}V @VVdV \|\\
 0 @>>> C@>>i> D @>>> Z@>>>0\equiv G\\
 & &@VVV @VVV\\
 & &0 &&0\\
&& |||&&|||\\
 & &H && V.
\end{CD}
$$
\begin{enumerate}
\item[i)] $V$ and $G$ are $\mathcal A$-trivial if and only if $F$
and $H$ are $\mathcal A$-trivial. \item[ii)] If $F$ is $\mathcal
A$-trivial then $G$ is $\mathcal A$-trivial; if $V$ is $\mathcal
A$-trivial, $H$ is $\mathcal A$-trivial.
\end{enumerate}
\end{lemma}\begin{proof}The second part of assertion $ii)$ follows from $(1)$
in the previous Proposition \ref{pbct}. Let $A\in \mathcal A$ and
notice that from the diagram in the hypothesis it is immediate the
construction of the commutative diagram:
$$\begin{CD}
 & & & & &&\\
 & & & &@AAA @AAA\\
 & & & &\mathcal L(B,A)& = &\mathcal L(B, A)& \\
 & & & &@AA{b^*}A  @AA{a^*}A\\
 0@>>> \mathcal L(Z, A) @>>> \mathcal L(X, A) @>{j^*}>> \mathcal L(Y, A) @>>> \\
 & & \| & & @AA{d^*}A @AA{c^*}A \\
 0 @>>> \mathcal L(Z, A)@>>> \mathcal L(D, A) @>{i^*}>> \mathcal L(C, A)@>>>\\
 & & & &@AAA @AAA\\
 & & & &0 &&0.\\
\end{CD}$$\smallskip

\noindent Now, the first part of ii) can be easily obtained
chasing the diagram while i) follows simply observing that the
restriction operators $b^*$ and $i^*$ are surjective if and only
if $j^*$ and
$a^*$ are surjective. 
\end{proof}

\noindent \textbf{Remark.} Pe\l czy\'{n}ski's Proposition 2.6 of
\cite{pelcaverage} can be considered a rudimentary version of this
principle.\medskip

Returning to more classical terms, it is especially interesting
for us the characterization of $C(K)$-trivial extensions that
Zippin formulates and proves in \cite{zippill}.

\begin{lemma}\label{zips} A sequence $0\to Y\stackrel{j}\to X\to Z\to
Z\to 0$ is $(\lambda, C(K))$-trivial if and only if there is a
w*-continuous map $\omega: B_{Y^*} \to \lambda B_{X^*}$ such that
$j^* \omega = id$.
\end{lemma}

The map $\omega$ shall be called a $\lambda$-$w^*$-selector for
$j^*$. Zippin \cite{zippln,zipproc} uses this criterion to obtain
different proofs of the Lindenstrauss-Pe{\l}czy\'nski theorem. It
is inspired by the most natural possible situation: the sequence $
0 \to Y \stackrel{\delta_Y}\to C(B_{Y^*}) \to C(B_{Y^*})/Y \to 0
\equiv \complement_Y$ , in which $\delta_Y: Y \to C(B_{Y^*})$ is
the canonical embedding, is $(1, C(K))$-trivial; indeed, the map
$\omega:B_{Y^*}\to B_{C(B_{Y^*})}$ defined as
$\omega(x^*)(f)=f(x^*)$ is a $1$-$w^*$-selector for $\delta_Y^*$.
Let us remark that every $C(K)$-trivial sequence is a pull-back of
$\complement_Y$ and conversely.\medskip

It shall be useful to notice that some properties of the
w*-topology in $l_p$, $1\le p \leq \infty$, pass to $l_p$-vector
sums of Banach spaces. Given an $l_p$-sum $l_p(X_n)$ we denote
$\pi_j:l_p(X_n)\to X_j$ the natural projections.\smallskip

\begin{lemma}\label{wstar}
Let $(E_n^*)_n$ be a sequence of dual spaces and let $(x_k)_k$ be
a bounded sequence in $l_p(E_n^*)$, $1\le p \leq \infty$. The
sequence $(x_k)_k$ is $w^*$-null if and only if the sequences
$(\pi_j(x_k))_k$ are $w^*$-null.\end{lemma}\smallskip

\begin{proof} It is clear that if $(x_k)$ is $w^*$-null, the sequences $(\pi_j(x_k))_k$ are also
$w^*$-null. Conversely, let $x$ be an element of $l_{p^*}(E_n^*)$
--we understand $l_{1}(E_n^*)$ (resp. $c_0(E_n^*)$) when
$p=\infty$ (resp. $p=1$)-- that we write as $x= \lim s_n$ in such
a way that $(s_n)$ are finitely supported (that is, $\pi_j(s_n)=0$
except for a finite quantity of indices $j$). If for all $j$ the
sequence $(\pi_j(x_k))_k$ is $w^*$-null in $E_j$, fixed
$\varepsilon/2>0$, there exist $N_1>0$ and $N_2>0$ such that for
all $n>N_1$ one has $|x_k(x-s_n)|\le \varepsilon/2$, and for all
$k>N_2$ one has $|x_k(s_n)|\le \varepsilon/2$. Thus, $|x_k(p)|\le
\varepsilon$.\end{proof}

Given a family of exact sequences $0\to A_n\to B_n\to C_n\to
0\equiv F_n$, we call $0\to l_p(A_n)\to l_p(B_n)\to l_p(C_n)\to
0\equiv l_p(F_n)$ the $l_p$-amalgam of $(F_n)_n$, for $1 \leq
p\leq \infty$. Analogously, the $c_0$-amalgam of $(F_n)$ shall be
denoted $c_0(F_n)$. One has:

\begin{proposition}\label{amalgam}
The $c_0$- and $l_p$-amalgams of $(\lambda, C(K))$-trivial exact
sequences, $1\leq p < \infty$, are $(\lambda, C(K))$-trivial.
\end{proposition}

\begin{proof} For each $n$, let $\omega_n: B_{A_n^*} \to\lambda B_{B_n^*}$ be a
$\lambda$-$w*$-selector for $j_n^*$. If we have the
$l_p$-amalgam$$\begin{CD}0@>>> l_p(A_n)@>\chi>> l_p(B_n)@>>>
l_p(C_n)@>>>0\equiv l_p(F_n),\end{CD}$$ it follows from Lemma
\ref{wstar} that the map $\Omega:B_{l_{p^*}(A_n^*)} \to \lambda
B_{l_{p^*}(B_n^*)}$ defined by $\Omega[(a_n^*)]
=[\omega_n(a_n^*)]$ is a $\lambda$-$w^*$-selector for $\chi^*$.
\end{proof}

The situation for $l_\infty$-amalgams is entirely different
because a subspace $X$ of $l_\infty$ can only be
$C(K)$-complemented if it enjoys the Dunford-Pettis and the
Grothendieck character of $l_\infty$: indeed, every operator
defined on $X$ with separable range must be weakly compact, and
thus every weakly compact operator $X \to c_0$ must be completely
continuous. This means that, subspaces such as $l_\infty(l_2^n)$
cannot be $C(K)$-complemented in $l_\infty$ and therefore the
$l_\infty$-amalgam of the sequences

$$\begin{CD} 0@>>> l_2^n @>>> l_\infty^{m(n)} @>>>
l_\infty^{m(n)}/l_2^n @>>> 0 \end{CD}$$

in which the embeddings are $(1+\varepsilon)$-isometries cannot be
$C(K)$-trivial. Those sequences are
$((1+\varepsilon)^2,C(K))$-trivial by the existence of the
Bartle-Graves continuous selection.

Later on we shall also show that the $c_0$ -amalgam of $\mathcal
L_\infty$-trivial sequences is not necessarily $\mathcal
L_\infty$-trivial. Our concern now is to study the stability of
$C(K)$-trivial sequences by duality. Let us start observing that,
for every subspace $H$ of $c_0$ and every separably Banach space
$S$, the sequence $0\to H\to X\to S\to 0\equiv F$ is
$C(K)$-trivial as well as its dual $F^*$. That $F$ is
$C(K)$-trivial was observed by Lindenstrauss and Pe\l czy\'{n}ski
in \cite[Cor. 2]{lindpelc}. The assertion about $F^*$ (actually,
that every exact sequence $0\to Y \to X \to H^* \to 0$ is
$C(K)$-trivial for every subspace $H$ of $c_0$) directly follows
from the equality $\Ext(H^*, C(K))=0$, which is implied by the
Johnson-Zippin theorem. It is not true, in general, that the dual
or bidual of a $C(K)$-trivial sequence is $C(K)$-trivial:\\

\noindent \textbf{Examples.} Consider the $C(K)$-trivial sequence
$ 0 \to l_2 \stackrel{\delta}\to C(B_{l_2})  \to Q \to 0 \equiv F
$. We claim that the dual sequence $0\to  Q^* \to L_1
\stackrel{\delta^*}\to l_2 \to 0 \equiv F^* $ cannot be
$C(K)$-trivial. So, assume it to be otherwise. Consider a
projective presentation of $l_2$ $$
\begin{CD}
0@>>> K(l_2) @>>> l_1  @>>> l_2 @>>>0 \equiv P.
\end{CD}
$$The space $K(l_2)$ is complemented in its bidual as it was
proved by Kalton and Pe\l czy\'{n}ski in \cite{kaltpeck}. Hence
$\Ext(L_1, K(l_2))=0$ using Lindenstrauss's lifting principle.
Thus, the lower pull-back sequence in the diagram $$
\begin{CD}
0@>>> K(l_2) @>>> l_1  @>>> l_2 @>>>0 \equiv P\\ &&\Vert && @AAA
@AA{\delta^*}A\\ 0 @>>> K(l_2) @>>> PB  @>>> L_1 @>>>0 \equiv P
\delta^*
\end{CD}
$$splits, and therefore the quotient map $L_1 \to l_2$ can be
lifted to $l_1$. This lifting can be chosen so that its
restriction to $Q^*$ gives a surjective operator $\phi: Q^* \to
K(l_2)$. In this way $P\equiv \phi F^*$ (see also
\cite{castmoreisr}) and $P$ a push-out of $F^*$. So, it must be
$C(K)$-trivial. If a projective presentation $P$ of $l_2$ is
$C(K)$-trivial then $\Ext(l_2, C(K))=0$. But it was proved by
Kalton in \cite{kaltc0}, see also \cite{ccky}, that $\Ext(l_2,
C[0,1])\neq 0$. The bidual sequence $$
\begin{CD}
0@>>> l_2 @>>> I @>>> Q^{**} @>>>0 \equiv F^{**}
\end{CD}
$$cannot be $C(K)$-trivial neither: since $I$ is injective, every
operator $I \to C(B_{l_2})$ is weakly compact, hence completely
continuous; thus, the canonical inclusion $l_2 \to C(B_{l_2})$
cannot be extended to $I$. Another example is provided by
sequences having the form$$
\begin{CD}
0@>>> c_0(A_n) @>>> c_0(l_\infty^{m(n)})   @>>> c_0 (C_n) @>>>0
\equiv F;
\end{CD}
$$they satisfy that $F$ and $F^*$ do $C(K)$-split although
$F^{**}$ does not necessarily $C(K)$-split.

\section{On Banach spaces of Lindenstrauss-Pe\l czy\'{n}ski type}

There is an obvious difference between the
Lindenstrauss-Pe{\l}czy\'{n}ski and the Johnson-Zippin theorems.
While the former asserts that every sequence $0\to H\to c_0\to
c_0/H\to 0\equiv F$ is $C(K)$-trivial the latter establishes that
the dual sequence $F^*$ is $\mathcal L_\infty$-trivial. Let us see
that the class $\mathcal L_\infty$ cannot be enlarged, obtaining
in this form a new characterization of $\mathcal L_\infty$-spaces.

\begin{proposition}\label{jz} For a Banach space $E$ the following are equivalent:
\begin{enumerate}
\item $E$ is an $\mathcal L_\infty$-space. \item Every $E$-valued
operator defined on a $w(l_1,c_0)$-closed subspace of $l_1$ can be
extended to $l_1$.
\end{enumerate}
\end{proposition}
\begin{proof}
It is well-known that (2) can be written as: $\Ext(H^*,E)=0$ for
every subspace $H$ of $c_0$. Using again the decomposition $0\to
c_0(A_n) \to H \to c_0(B_n) \to 0$, where $A_n, B_n$ are
finite-dimensional spaces, plus a simple 3-space argument (see
\cite[Cor. 1.2]{ccky}, the previous condition is equivalent to:
$\Ext(l_1(G_n),E)=0$ for every sequence $(G_n)$ of finite
dimensional spaces. It was already observed by Johnson \cite{john}
that a sequence $0 \to Y \to X \to Z \to 0 \equiv F$ locally
splits if and only if $FT\equiv 0$ for every operator $T: l_1(G_n)
\to Z$, with $G_n$ finite-dimensional. So, every sequence $0 \to E
\to X \to Z \to 0$ locally splits and in particular it does any
sequence $0 \to E \to l_\infty(I) \to Q \to 0 \equiv G$. Recall
from \cite{kaltloc} that an exact sequence $F$ locally splits if
and only if $F^{**}$ splits. Hence, the bidual sequence $G^{**}$
splits, $E^{**}$ must be complemented in an $\mathcal
L_\infty$-space and $E$ must be itself an $\mathcal
L_\infty$-space.
\end{proof}

\enlargethispage*{50pt}Not entirely trivial is the observation
that condition (2) can be replaced by (2') For every set
$\Gamma)$, every $E$-valued operator defined on a
$w(l_1(\Gamma),c_0(\Gamma))$-closed subspace of $l_1(\Gamma)$ can
be extended to $l_1(\Gamma)$. The proof only requires to take in
consideration the decomposition lemma of \cite{johnzipproc}.

The situation outlined for $w^*(l_1,c_0)$-closed subspaces of
$l_1$, together with the Lindenstrauss-Pe\l czy\'{n}ski's remark
in \cite{lindpelc} asserting that operators with range isometric
preduals of $L_1$ extend from subspaces of $c_0$ to the whole
space, suggest to investigate how much the class of $C(K)$-spaces
can be enlarged in the Lindenstrauss-Pe\l czy\'{n}ski
theorem.\smallskip

\begin{defin}[Spaces of Lindenstrauss-Pe{\l}czy\'{n}ski type] We shall
say that a Banach space $E$ is a Lindenstrauss-Pe{\l}czy\'{n}ski
space, in short an $\mathscr {LP}$ space, if all operators from
subspaces of $c_0$ into $E$ can be extended to $c_0$.\end{defin}
We shall also need the quantitative version: when every operator
$T:H\to E$ admits an extension $\widehat{T}:c_0\to E$ such that
$\|\widehat{T}\|\leq \lambda \|T\|$ we shall say that $E$ is an
$\mathscr{LP}_\lambda$ space. It is not hard to see that every
$\mathscr{LP}$ space is an $\mathscr{LP}_\lambda$ space for some
$\lambda$.

\begin{lemma} Every $\mathscr{LP}_\lambda$ space is an $\mathcal
L_{\infty, 2\lambda}$ space.\end{lemma}
\begin{proof} Let $E$ be an $\mathscr{LP}_\lambda$ space.
Let $T: Y \to E$ be a compact operator from a subspace $Y$ of a
separable space $X$. Then $T$ factorizes as through some subspace
$i:H \to c_0$ as $T=BA$ with $A: Y \to H$ and $B: H \to E$. By
definition, there is an extension $B_1: c_0 \to E$ of $B$ with
$\|B_1\|\leq \lambda \|B\|$; while Sobczyk's theorem gives an
extension $A_1: X \to c_0$ of $iA$ with $\|A_1\|\leq 2\|iA\|$. The
composition $B_1A_1: X \to E$ extends $T$ and verifies
$\|B_1A_1\|\leq 2\lambda\|T\|$. Using Lindenstrauss's
characterization \cite{lindmem}, $E$ must be an $\mathcal
L_{\infty,2\lambda}$ space.\end{proof}

The converse fails: we show that not every $\mathcal L_\infty$
space is an $\mathscr{L P}$ space. This solves Zippin's Problem
6.15 in \cite{zipphandbook}. The example (which was sketched in
\cite{ccmy}) is based on the Bourgain-Pisier construction
\cite{bourpisi} which shows that for every separable Banach space
$X$ there is an exact sequence
$$\begin{CD}0@>>> X @>>> \mathcal L_\infty(X)@>>>\mathcal
L_\infty(X)/X @>>> 0,
\end{CD}$$
in which $\mathcal L_\infty(X)$ is a separable $\mathcal
L_\infty$-space and $\mathcal L_\infty(X)/X$ has the Schur
property.
\begin{proposition} Let $H$ be a subspace of $c_0$ such that $c_0/H$ is not isomorphic to $c_0$. Then
$\mathcal L_\infty(H)$ is not an $\mathscr{LP}$
space.\end{proposition}
\begin{proof} Let us consider a sequence $0 \to H \stackrel{j}\to c_0 \to c_0/H \to 0\equiv F$. Let $0\to
H\stackrel{i}\to \mathcal L_\infty(H)\to S\to 0\equiv \mathcal
H_\infty$ be the Bourgain-Pisier sequence. Assume that $i$ extends
to $c_0$ through $j$, so that $F$ is a pull-back of $\mathcal
H_\infty$. By Sobczyk's theorem, $\mathcal H_\infty$ is a
pull-back of $F$. Applying the first diagonal principle of
\cite{lindrose} one gets an isomorphism$$\mathcal L_\infty(H)
\oplus c_0/H \simeq c_0 \oplus S$$In particular, $c_0/H$ is a
complemented subspace of $c_0\oplus S$. Since $S$ and $c_0$ are
totally incomparable by the Schur property of $S$ the
decomposition theorem of Edelstein-Wojtasczyk again (see
\cite[Theorem 2.c.13]{lindtzaf}) ensures that $c_0/H$ is
isomorphic to some $A\oplus B$ with $A$ complemented in  $c_0$ and
$B$ complemented in $S$. Since $c_0/H$ is a subspace of $c_0$, $B$
can only be finite dimensional, hence $c_0/H\backsimeq c_0$,
against the hypothesis.\end{proof}

This example immediately implies

\noindent \textbf{Example.} \emph{The $c_0$-amalgam of $(\mu,
\mathcal L_\infty)$-trivial sequences is not necessarily $\mathcal
L_\infty$-trivial.}
\begin{proof} General structure results of Johnson-Rosenthal and Zippin (see \cite[1.g.2 and
2.d.1]{lindtzaf}) imply that given a subspace $H$ of $c_0$ there
exist sequences $(A_n)$ and $(B_n)$ of finite dimensional spaces
such that there is an exact sequence $ 0 \to c_0(A_n) \to H \to
c_0(B_n) \to 0$. There is little doubt that exact sequences $ 0\to
A_n  \to l_\infty^{m(n)} \to C_n \to 0$ are $(\lambda, \mathcal
L_{\infty, \lambda})$-trivial. If all the amalgams $ 0\to c_0(A_n)
\to c_0(l_\infty^{m(n)}) \to c_0(C_n) \to 0$ were $\mathcal
L_\infty$-trivial then every sequence $ 0\to H\to c_0 \to c_0/H
\to 0$ should also be $\mathcal L_\infty$-trivial: indeed, there
would be a complete push-out diagram$$
\begin{CD}
 & &0 &&0\\
 & &@VVV @VVV\\
 & &c_0(A_n)&=& c_0(A_n)\\
 & &@VV{a}V  @VV{i}V\\
 0@>>> H @>j>>c_0 @>>> c_0/H @>>> 0\equiv F\\
 & &@VV{b}V @VVcV \|\\
 0 @>>> c_0(B_n)@>>d> Z @>>> Z/c_0(B_n) @>>>0\equiv bF\\
 & &@VVV @VVV\\
 & &0 &&0\\
&& |||&&|||\\
 & &Vd && V,\\
\end{CD}
$$to which we apply i) of Lemma \ref{cruz}: if $V$ and $bF$ are
$\mathcal L_\infty$-trivial then $F$ is $\mathcal
L_\infty$-trivial (which we know it is not). It remains to check
that $V$ and $bF$ are $\mathcal L_\infty$-trivial. The sequence
$V$ is $\mathcal L_\infty$-trivial by our assumption and the
Lindenstrauss-Rosenthal theorem that asserts that $c_0$ is
automorphic (see \cite{lindrose,castmoreisr}); i.e., there exists
an isomorphism $\tau: c_0 \to c_0$ making a commutative diagram
$$
\begin{CD}
 0 @>>> c_0(A_n)@>>> c_0 @>>> Z @>>>0\\
 & &\| &&@VV{\tau}V  @VVV\\
 0 @>>> c_0(B_n) @>>> c_0(l_\infty^{m(n)}) @>>>  c_0(l_\infty^{m(n)}/A_n) @>>>0\\
\end{CD}
$$and therefore the two sequences $\mathcal L_\infty$-split
simultaneously. The sequence $bF$ is $\mathcal L_\infty$-trivial
with essentially the same arguments taking into account that $Z$
must be a subspace of $c_0$.\end{proof}

The problem of identifying $\mathscr{LP}$ spaces is still far from
being solved, and it actually gives rise to interesting questions.
Observe that, in addition to $C(K)$-spaces, it is clear that
complemented subspaces of $C(K)$-spaces and separably injective
spaces are also $\mathscr{LP}$ spaces. The reader might be
surprised by the distinction between the two, especially regarding
the fact that every injective space is complemented in some
$C(K)$-space. Let us show that the distinction is necessary.

\begin{proposition}
There exists a separably injective space that is not complemented
in any $C(K)$-space.\end{proposition}\begin{proof} Let us consider
the pull-back diagram$$
\begin{CD}
0@>>> c_0 @>>>l_\infty @>>> l_\infty/c_0 @>>>0\equiv \mathcal I\\
&&\parallel &&@AAA @AA{\lambda \cdot 1}A \\ 0@>>> c_0 @>>>
P(\lambda) @>>> l_\infty/c_0 @>>>0\equiv \mathcal I_\lambda.
\end{CD}
$$Benyamini shows in \cite{beny} that $P(\lambda)$ is not less than
$\lambda^{-1}$-complemented in any $C(K)$-space. Thus, the
$c_0$-amalgam of the family  $(\mathcal I_{n^{-1}})$ $$
\begin{CD}
0@>>> c_0(c_0)@>>> c_0(P(n^{-1})) @>>> c_0(l_\infty/c_0)@>>>
0\equiv c_0(\mathcal I_{n^{-1}})
\end{CD}
$$ provides an exact sequence in which  both $c_0(c_0)$ as well as
$c_0(l_\infty/c_0)$ are $C(K)$-spaces. However, the space
$c_0(P(n^{-1}))$ cannot be complemented in any $C(K)$-space. That
$c_0$ is separably injective is precisely Sobczyk's theorem. That
$l_\infty/c_0$ is separably injective is well known and follows
from Proposition \ref{sepinj} below. It was shown in
\cite{johnoikh,rose} that when $X$ is separably injective then
$c_0(X)$ is separably injective as well. Finally, separable
injectivity is a 3-space property (see \cite[Cor.
1.2]{ccky}).\end{proof}

There are other $\mathscr{LP}$ spaces. As we mentioned before,
according to \cite[remark 2, p.234]{lindpelc}, isometric preduals
of $L_1$ are $\mathscr{LP}$ spaces. So, it is quite natural to ask
whether the previous classes (namely: complemented subspaces of a
$ C(K)$-space, separably injective spaces and isometric preduals
of $L_1$) exhaust the $\mathscr{LP}$ spaces. The answer is no.

\begin{proposition}\label{tipo0} Every $\mathcal L_\infty$-space not containing $c_0$
is an $\mathscr{LP}$ space.
\end{proposition}
\begin{proof} Let $X$ be a Banach space. It is not difficult to
see that when $Y$ contains no copy of $c_0$ then every operator $X
\to Y$ is unconditionally converging (see \cite{djt}). As it is
well known, every subspace $H$ of $c_0$ has Pe{\l}czy\'nski
property $(V)$, which means that every unconditionally converging
operator $H \to X$ is weakly compact. Since $H^*$ is Schur, every
weakly compact operator $H \to X$ must be compact; hence, when $Y$
contains no copy of $c_0$ every operator $H \to Y$ must be
compact. Using Lindenstrauss's extension theorem for compact
operators \cite{lindmem}, the result follows.\end{proof}

By a result of Johnson and Zippin \cite{johnzippre} separable
isometric $L_1$-preduals are quotients of $C[0,1]$. Observe that
the $\mathcal L_\infty$-spaces not containing $c_0$ cannot be
quotients of $C(K)$-spaces. Concrete examples of $\mathcal
L_\infty$-spaces  not containing $c_0$ can be obtained applying
the Bourgain-Pisier construction $0\to X \to \mathcal
L_\infty(X)\to S \to 0$ to spaces $X$
 without copies of $c_0$ by a simple 3-space argument (see
\cite{castgonz}).\medskip

Recall (see \cite{castgonz}) that a property $P$ is said to be a
3-space property if whenever the spaces $Y$ and $Z$ in an exact
sequence \seq have $P$ then also $X$ has $P$. New classes of
$\mathscr{LP}$ spaces can be obtained showing that this class
satisfies the 3-space property.

\section{The 3-space problem for $\mathscr{LP}$ spaces}

The purpose of this section is to show:

\begin{theorem}\label{3sp} The class of $\mathscr{LP}$-spaces has the
3-space property.
\end{theorem}

The proof however is not simple and requires both a different
characterization of $\mathscr{LP}$-spaces and a new method to
obtain 3-space properties. We assume from the reader some a
cquaintance with the theory of operator ideals such as developed
by Pietsch in \cite{piet}. Recall that an operator ideal
$\mathfrak A$ is said to be surjective (see \cite[4.7.9]{piet}) if
whenever $Q$ is a quotient map and $TQ \in \mathfrak A$ imply $T
\in \mathfrak A$; dually, the ideal $\mathfrak A$ is injective
(see \cite[4.6.9]{piet}) if whenever $J$ is an into isomorphism
and $JT\in \mathfrak A$ imply $T \in \mathfrak A$. Let us consider
the operator ideal $\mathfrak{J}_{0}$ of those operators that
factorize through a subspace of $c_0$.

\begin{prop}
A Banach space $E$ is an $\mathscr{LP}$ space if and only if the
functor $\mathfrak{J}_{0}(\cdot,E)$ is exact when acting on the
category of separable Banach spaces.
\end{prop}
\begin{proof} Let $j:Y \to X$ be an into isomorphism, and let $T \in \mathfrak J_0(Y, E)$. Assume that
$T=RS$ with $S\in \mathfrak L(Y,H)$, $R \in \mathfrak L(H, E)$ and
$i: H \to c_{0}$ an into isomorphism. The operator $R$ can be
extended to an operator $R_{1}:c_{0}\to E$ through $i$ since $E$
is an $\mathscr {LP}$-space; while Sobczyk's theorem allows one to
extend $iS: Y \to c_0$ to an operator $S_1: X \to c_0$ through
$j$. All together provides an extension $R_1 S_1$ of $T$ through
$j$.  The other implication is immediate.
\end{proof}

The new method to obtain 3-space properties is the following.

\begin{proposition}\label{metodo} Let $\mathfrak A$ be a surjective and injective operator ideal. The
class of all Banach spaces $E$ such that the functor $\mathfrak
A(\cdot, E)$ is exact has the 3-space property.
\end{proposition}
\begin{proof} The surjectivity of $\mathfrak{A}$ yields that given an exact sequence $0 \to Y \stackrel{j}\to X \stackrel{q}\to Z \to 0$
and a given space $E$  the induced sequence$$
\begin{CD} 0 @>>>\mathfrak{A}(Z,E)@>q^*>>\mathfrak{A}(X,E)@>j^*>>\mathfrak{A}(Y,E)\\
\end{CD}$$is exact. Now, let$$\begin{CD} 0 @>>>A@>i>>B@>p>>C@>>>0\\
\end{CD}$$be an exact sequence. By assumption, both $\mathfrak{A}(\cdot,A)$ and $\mathfrak{A}(\cdot,C)$ are exact functors
and we need to prove that also $\mathfrak{A}(\cdot,B)$ is exact.
To this end we construct the commutative diagram

$$\begin{CD}
 && 0&&0&&0 \\
 && @VVV @VVV @VVV \\
0@>>>\mathfrak{A}(Z,A)@>q^{*}>>\mathfrak{A}(X,A)@>j^{*}>>\mathfrak{A}(Y,A)@>>>\mathfrak{A}(Y,A)\diagup
j^{*}( \mathfrak{A}(X,A))\\ && @Vi^*VV @Vi^*VV @Vi^*VV @VVV\\ 0
@>>>\mathfrak{A}(Z,B)@>q^{*}>>\mathfrak{A}(X,B)@>j^{*}>>\mathfrak{A}(Y,B)@>>>\mathfrak{A}(Y,B)\diagup
j^{*}(\mathfrak{A}(X,B))\\ && @Vp^*VV  @Vp^*VV  @Vp^*VV @VVV\\ 0
@>>>\mathfrak{A}(Z,C)@>q^{*}>>\mathfrak{A}(X,C)@>j^{*}>>\mathfrak{A}(Y,C)@>>>\mathfrak{A}(Y,C)\diagup
j^{*}(\mathfrak{A}(X,C)).
\end{CD}$$
The rows are exact by the surjectivity of $\mathfrak U$, while the
columns are also exact by injectivity of $\mathfrak U$. By
hypothesis, $$\mathfrak{A}(Y,A)\diagup j^{*}( \mathfrak{A}(X,A)) =
 \mathfrak{A}(Y,C)\diagup j^{*}(
 \mathfrak{A}(X,C)) = 0$$ and the exactness of the fourth column implies that $$\mathfrak{A}(Y,B)\diagup j^{*}(
 \mathfrak{A}(X,B)) = 0,$$ hence $\mathfrak{A}(\cdot,B)$ is exact.
\end{proof}

Since the separability assumption of the previous characterization
of $\mathscr {LP}$-spaces does not affect the method of
Proposition \ref{metodo}, the proof of Theorem \ref{3sp} will be
complete after showing:

\begin{lema} The ideal $\mathfrak{J}_{0}$ is injective and surjective.
\end{lema}
\begin{proof} The injectivity is a direct consequence of the
definition. To show the surjectivity, let $\tau: X \to E$ be an
operator which factorizes as $\tau = \varphi_0\varphi_1$ through a
subspace $H$ of $c_0$ in a diagram $$\begin{CD} 0
@>>>Y@>j>>X@>p>>Z@>>>0\\
 &&&& @V\varphi_{1}VV\\
 &&&& H\\
 &&&& @V\varphi_{0}VV\\
 &&&& E
 \end{CD}$$
Assume that $\tau j=0$. One then has the commutative diagram $$
\begin{CD}
0 @>>>Y@>j>>X@>p>>Z@>>>0\\
 &&@VVV  @V\varphi_{1}VV @V\widetilde{\varphi}_{1}VV\\
0@>>> \ker \varphi_{0}@>>>H@>>P> H / \ker \varphi_{0} @>>> 0\\
 &&&& @V\varphi_{0}VV\\
 &&&& E
 \end{CD}$$\\

It is clear that there exists an operator
$\widetilde{\varphi}_{0}:H / \ker \varphi_{0} \rightarrow E$ such
that $\widetilde{\varphi}_{0}P= \varphi_0$. It is then obvious
that  $\widetilde{\varphi}_{0}\widetilde{\varphi}_{1} p
={\varphi}_{0}{\varphi}_{1}$. Moreover the operator
$\widetilde{\varphi}_{0}\widetilde{\varphi}_{1}\in
\mathfrak{I}_{0}(Z,E)$ since $H/\ker \varphi_{0}$, as a quotient
of a subspace of $c_{0}$, is itself a subspace of $c_{0}$ (see
\cite{lindtzaf}).
\end{proof}

Theorem \ref{3sp}, in particular, yields:

\begin{cor} Every twisted sum of $C(K)$-spaces is an $\mathscr
{LP}$-space.
\end{cor}

The paper \cite[Thm. 2.3]{ccky} contains most of the available
information about how to construct twisted sums of $C(K)$-spaces.
For instance, it is shown that for every separable Banach space
$X$ not containing $l_1$ there exists an exact sequence $$
\begin{CD}
0@>>> C[0,1] @>i>> \Omega(X) @>q>> X@>>> 0\end{CD}$$ with strictly
singular quotient map. Of course, the space $\Omega(X)$ is not a
quotient of a $C(K)$-space. Using Theorem 4.7 in \cite{ccky} one
can obtain examples of $\mathscr{LP}$-spaces not containing $l_1$
which are not $C(K)$-spaces. We do not know if there exist
$\mathcal L_\infty$-spaces not containing $l_1$ which are not
$\mathscr{LP}$ spaces. See also Question 1 in Section 4.\\

Moving to the nonseparable spaces, new $\mathscr{LP}$ spaces can
be obtained through the following stability result.

\begin{proposition}\label{sepinj} Every quotient of an $\mathscr{LP}$ space by a
separably injective space is an $\mathscr{LP}$ space.
\end{proposition}\begin{proof}Let us consider an exact sequence
$ 0 \to SI\to \mathscr {LP} \stackrel{q}\to Q \to 0$ in which the
middle space is a Lindenstrauss-Pe\l czy\'{n}ski space and the
subspace is separably injective. Let $\phi: H \to Q$ be an
operator from a subspace $H$ of $c_0$. Since $SI$ is separably
injective, $\Ext(H, SI)=0$. Hence $F\phi$ splits and $\phi$ can be
lifted through $q$ to an operator $\psi: H \to \mathscr{LP}$. This
operator can be extended to an operator $\Psi: c_0 \to
\mathscr{LP}$. The operator $q \Psi: c_0 \to Q$ is the desired
extension of $\phi$.\end{proof}The same proof provides that the
quotient of two separably injective spaces is separably injective.
In particular, $l_\infty/c_0$ is separably injective. We would
like to mention one more stability result. Being obvious that
$l_\infty$-vector sums of sequences of $\mathscr{LP}_\lambda$
spaces are $\mathscr{LP}$ spaces, the corresponding result for
$c_0$-vector sums keeps being true. Details shall appear
elsewhere.

\begin{proposition}
Let $\lambda\geq 1$. The $c_0$-vector sum of a sequence of
$\mathscr{LP}_\lambda$ spaces is an $\mathscr{LP}$ space.
\end{proposition}

\section{Some open questions}

Regarding the stability result in Proposition \ref{sepinj} it
seems quite natural to ask:\smallskip

\noindent \textbf{Question 1.} Is the quotient of two $\mathscr
{LP}$ spaces an $\mathscr{LP}$ space?\smallskip

An affirmative answer to Question 1 would imply that separable
$\mathcal L_\infty$-spaces not containing $l_1$ are $\mathscr{LP}$
spaces. Two particularly interesting quotients of two
$\mathscr{LP}$-spaces are remarked in the next question.
\smallskip

\noindent \textbf{Question 2.} Is $l_\infty/C[0,1]$ an $\mathscr
{LP}$ space?  Must $\mathcal L_\infty$-spaces which are quotients
of $C[0,1]$ be $\mathscr{LP}$ spaces ?\medskip

As we have already mentioned in the introduction, Johnson and
Zippin proved in \cite{johnzipproc} that every extension $0\to
H\to c_0(\Gamma)\to Z\to 0$ is $C(K)$-trivial. It will be nice to
know if $\mathscr{LP}$ spaces can play the role of $C(K)$-spaces
in this result.\smallskip

\noindent \textbf{Question 3.} Given a subspace $H$ of
$c_0(\Gamma)$, does every operator $H\to \mathscr{LP}$ have an
extension to $c_0(\Gamma)$?\smallskip

Needless to say, the extension property one would like to get from
$\mathscr{LP}$ spaces is: every $C(K)$-trivial sequence is also
$\mathscr{LP}$-trivial. Unfortunately, this does not
hold.\smallskip

\noindent \textbf{Example}. We already know that the
Bourgain-Pisier space $\mathcal L_\infty(l_2)$ is an
$\mathscr{LP}$ space. Consider then the sequences $$
\begin{CD}
0 @>>> l_2 @>{\delta_2}>> C(B_{l_2}) @>>> Q @>>> 0 \\
 && \parallel &&& &\\
0@>>>l_2 @>>j> \mathcal L_\infty(l_2) @>>> S @>>> 0. \end{CD} $$If
$j$ could be extended to an operator $J: C(B_{l_2}) \to \mathcal
L_\infty(l_2)$ through $\delta_2$ this would be a weakly compact
operator since $\mathcal L_\infty(l_2)$ does not contain $c_0$.
Hence $J$ would be completely continuous by the Dunford-Pettis
property of $C(K)$-spaces. It is therefore impossible that
$J\delta_2 =j$.\medskip

The method of proof developed in Proposition \ref{metodo} is new.
It is moreover clear that it can be applied to other injective and
surjective operator ideals appearing in the literature. The
following ideals are injective and surjective (see \cite{piet})
$\mathfrak{L}$ = all operators; $\mathfrak{F}$ =finite rank
operators; $\mathfrak{K}$ = compact operators; $\mathfrak{W}$ =
weakly compact operators; $\mathfrak{U}$ = unconditionally summing
operators; $\mathfrak{J}_2$ = operators factorable through a
Hilbert space. Let us introduce some notation: given an injective
and surjective operator ideal $\mathfrak U$ let $E(\mathfrak U)$
be the class of all Banach spaces $E$ such that the functor
$\mathfrak U(\cdot, E)$. From Proposition \ref{metodo} we have
obtained easy proofs that the following classes have the 3-space
property: $E(\mathfrak L)$ = injective spaces; applying the method
only to separable spaces one obtains the class of separably
injective spaces; $E(\mathfrak K) = \mathcal L_\infty$-spaces (by
\cite{lindmem}); $E(\mathfrak W)$ = $\mathfrak L_\infty$-spaces
with the Schur property (shown in \cite{bourdelb}). The classes
$E(\mathfrak U)$ and $E(\mathfrak J_2)$ seem not have been
characterized yet.

A simple homological duality argument yields:

\begin{proposition} Let $\mathfrak A$ be a injective and surjective
operator ideal. The class of all Banach spaces $E$ such that the
functor $\mathfrak A(E, \cdot)$ is exact has the 3-space property.
\end{proposition}

If we call $\exists(\mathfrak{A})$ the previous class determined
by the ideal $\mathfrak A$ then the only non-trivial case
identified is $\exists(\mathfrak{K})= \mathcal L_1$-spaces.

\end{document}